\documentclass[11pt]{amsart}

\usepackage{amssymb,latexsym,verbatim,pdfsync,color}

\def\bC{{\mathbf C}}

\def\bR{{\mathbf R}}

\def\cE{{\mathcal E}}
\def\cF{{\mathcal F}}
 
\def\cH{{\mathcal H}}

\def\cK{{\mathcal K}}

\def\cM{{\mathcal M}}

\def\cO{{\mathcal O}}

\def\cR{{\mathcal R}}

\def\cZ{{\mathcal Z}}

\def\Re{\operatorname{Re}}
\def\Im{\operatorname{Im}}

\def\la{\langle}
\def\ra{\rangle}

\def\eps{\varepsilon}
\def\z{\zeta} 
\def\vp{\varphi}
\def\ov{\overline}
\def\p{\partial}

\def\ms{\medskip}

\def\MS{M\"untz--Sz\'asz}
\def\tnt{\textstyle{\frac{n}{2}}}
\def\tkt{\textstyle{\frac{k}{2}}}

\def\Hol{\operatorname{Hol}}

\def\rz{\big(1/z_j\big)}
\def\rw{\big(1/w_j\big)}

\newtheorem{thm}{Theorem}[section]
\newtheorem{lem}[thm]{Lemma}
\newtheorem{prop}[thm]{Proposition}

\newtheorem{remark}[thm]{Remark}

\newtheorem{exam}[thm]{Example}

\newtheorem{theorem}{Theorem}[]
\newtheorem*{definition}{Definition}

\begin{document}

\title[Function
  spaces on a half-plane]{Functions of exponential growth in a
    half-plane, 
sets of uniqueness
and the M\"untz--Sz\'asz problem for the Bergman space}
\author[M. M. Peloso]{Marco M. Peloso}
\author[M. Salvatori]{Maura Salvatori}
\address{Dipartimento di Matematica, Universit\`a degli Studi di
  Milano, Via C. Saldini 50, 20133 Milano, Italy}
\email{{\tt marco.peloso@unimi.it}}
\email{{\tt maura.salvatori@unimi.it}}
\thanks{{\em Math Subject Classification} 	30H99, 46E22, 30C15, 30C40.}
\thanks{Authors partially supported by the grant Prin 2010-11 {\em
    Real and Complex Manifolds: Analysis, Geometry, Topology and Harmonic
    Analysis}.}
\begin{abstract}
We introduce and study some new spaces of
holomorphic  functions on
the right half-plane $\cR$.  In a previous work,  S. Krantz, C. Stoppato 
and the first named author
formulated the 
 \MS\ problem for the Bergman
space, that is, the problem to characterize the sets of complex powers
$\{\z^{\lambda_j-1}\}$ with $\Re\lambda_j>0$ that form a complete
set in the Bergman space $A^2(\Delta)$, where $\Delta=\{\z:\,
|\z-1|<1\}$.   

In this paper, we construct a space of holomorphic functions on
the right half-plane, that we denote
 by $\cM^2_\omega(\cR)$, whose sets of uniqueness $\{\lambda_j\}$ 
correspond exactly to the
sets of powers  $\{\z^{\lambda_j-1}\}$ that are a complete set in
$A^2(\Delta)$.

We show that $\cM^2_\omega(\cR)$ is a reproducing kernel Hilbert
space and we prove a Paley--Wiener type theorem among several other 
structural
properties.

We introduce a transform $M_\Delta$ modelled on the classical
Mellin transform and 
show that $\cM^2_\omega(\cR) = 
2^{-z}\Gamma(1+z) M_\Delta(A^2(\Delta))$.  We determine a sufficient
condition on a set $\{\lambda_j\}$ to be a set of uniqueness for
$\cM^2_\omega(\cR)$, thus providing a sufficient condition for the
solution of the 
\MS\ for the Bergman space.
\end{abstract}
\maketitle

\section*{Introduction and statement of the main results}

In this paper we introduce and begin the analysis of a space of
holomorphic  functions on
the right half-plane.  The initial motivation for such a study arose
in the work on Bergman spaces of worm domains in $\bC^2$ by  S. Krantz
and the first named author.  In collaboration also with C. Stoppato \cite{KPS2} 
we stated the 
 \MS\ problem for the Bergman
space 
 and proved a preliminary result.

We denote by $\Delta$  the disk $\{\z:\, |\z-1|<1\}$, by $dA$ the
Lebesgue measure in $\bC$ and consider the (unweighted) Bergman space 
$A^2(\Delta)$.  Then the complex powers $\{\zeta^{\lambda-1}\}$
with $\Re\lambda>0$ are well defined and in $A^2(\Delta)$.  We denote
by $\cR$ the right half-plane and by $\ov\cR$ its closure. 

Following \cite{KPS2}, the \MS\ problem for the Bergman
space is the question of characterizing the sequences $\{\lambda_j\}$ in
$\cR$ such that $\{\z^{\lambda_j-1}\}$ is a complete set in
$A^2(\Delta)$, that is, $\operatorname{span}\{\z^{\lambda_j-1}\}$
is dense in $A^2(\Delta)$.\ms

The classical \MS\ theorem concerns with the completeness of
a set of powers $\{t^{\lambda_j-\frac12}\}$ in $L^2\big([0,1]\big)$, 
where $\Re\lambda_j>0$. 
The solution was provided in two papers separate by C. M\"untz \cite{Muntz} and by
O. Sz\'asz \cite{Szasz} where they show that 
the set $\{t^{\lambda_j-\frac12}\}$  is complete $L^2\big([0,1]\big)$ if and only if the
sequence $\{\lambda_j\}$ is a set of uniqueness for the Hardy space 
of the right half-plane $H^2(\cR)$, that is, if $f\in H^2(\cR)$ and
$f(\lambda_j)=0$ for every $j$, then $f$ is identically $0$. 

As in the classical case,
in order to study the \MS\ problem for the Bergman space we wish to
transform the question into characterizing the sets of uniqueness for some
(Hilbert) space of holomorphic functions.  We now outline our
paradigm. 
\ms

For $f\in A^2(\Delta)$ and $z\in\cR$ we
define the Mellin--Bergman transform
\begin{equation}\label{M-Bt}
M_\Delta f(\lambda) = \frac{1}{\pi} \iint_\Delta f(\z) \ov \z^{\lambda-1}\, dA(\z)\, .
\end{equation}
The function $\z^{\ov \lambda-1}$ is well defined and belongs to
$A^2(\Delta)$.  Then a set  $\{\z^{\ov \lambda_j-1}\}$ is complete in
$A^2(\Delta)$ if and only if  $f\in A^2(\Delta)$ and
$M_\Delta f(\lambda_j)=0$ for all $j$ implies that $f$ vanishes identically.  

 Thus, the 
main task becomes to characterize the space 
$M_\Delta\big(A^2(\Delta)\big)$ and study its sets of uniqueness.
To this end we introduce a 
space of holomorphic functions on  $\cR$.   

\begin{definition}\label{def-cM2-omega}{\rm
For $0<b<\infty$, denote by $S_b$ the vertical strip
$\{z=x+iy: 0<x<b\}$ and by $H^2(S_b)$ the classical Hardy space
$$
H^2(S_b) =\big\{ f\text{ holomorphic in } S_b:\,
  \sup_{0<x<b} \int_{-\infty}^{+\infty} |f(x+iy)|^2\, dy
  <\infty\big\}\,.
$$

On $\ov\cR$ consider the Borel measure 
$\omega =\sum_{n=0}^{+\infty} \frac{2^n}{n!}\delta_{\frac{n}{2}}(x)\otimes
dy$, 
and the space $L^2(\ov\cR,d\omega)$;  explicitely, the 
norm is given by 
\begin{equation}\label{L2-omega-norm}
\|f\|_{L^2(\ov\cR,d\omega)}^2 =
\sum_{n=0}^{+\infty} \frac{2^n}{n!} \int_{-\infty}^{+\infty}
|f\big(\tnt +iy\big)|^2\, dy <\infty \,.
\end{equation}

We define $\cM^2_\omega(\cR)$ to be the space of
holomorphic functions $f$ on  $\cR$ such that:
\begin{itemize}
\item[(H)]
$f\in H^2(S_b)$ for every $0<b<\infty$;\smallskip
\item[(B)] $f\in L^2(\ov\cR,d\omega)$.
\end{itemize}
}
\end{definition}

Observe that
condition (H) implies
that $f$ admits a boundary value function in $L^2(\bR)$ and we take such
boundary values as the definition of $f$
 on the imaginary axis, which
is a positive  $\omega$-measure set.  More
formally, we could define $\cM^2_\omega(\cR)$ as the closure in
$L^2(\ov\cR,d\omega)$ of the functions satisfying (H) that admits
continuous extension to the closure of $\cR$.

We also point out that the measure $\omega$ has been found in a
constructive way and that it satisfies a uniqueness property, see
Theorem \ref{uniqueness-measure} below.
\ms

Observing that $\cM^2_\omega(\cR)$ is a subspace of $L^2(\ov\cR,d\omega)$, we prove
\begin{theorem}\label{cM2-Hilbert}
The space $\cM^2_\omega(\cR)$ is a Hilbert space with reproducing kernel
and with the unique inner product such that
$$
\|f\|_{\cM^2_\omega(\cR)}^2 =
\sum_{n=0}^{+\infty} \frac{2^n}{n!} \int_{-\infty}^{+\infty}
|f\big(\tnt +iy\big)|^2\, dy <\infty \,.
$$
\end{theorem}

Notice that trivially $H^2(\cR)\subset \cM^2_\omega(\cR)$ as a closed
subspace.  
Moreover, the
following simple facts hold true (see Proposition \ref{prop-1.2}):
\begin{itemize}
\item[(i)]
 the function $\Gamma(1+\delta z)\in \cM^2_\omega(\cR)$ for $0<\delta<1$,
  but $\Gamma(1+z)\not\in \cM^2_\omega(\cR)$;\smallskip
\item[(ii)] there exists $f$ holomorphic in $\ov\cR$ and satisfying 
(B), but $f\not\in
  H^2(S_b)$ if $b>\frac12$.
\end{itemize}

We study some basic structural properties of $\cM^2_\omega(\cR)$ and
we begin by proving the following Paley--Wiener-type theorem.  We
define 
the Fourier transform  of $\psi\in L^1(\bR)$
$$
(\cF \psi) (\xi) = \frac{1}{\sqrt{2\pi}} \int_{-\infty}^{+\infty}
\psi(x) e^{-ix\xi}\, dx\,. 
$$

\begin{theorem}\label{PW-thm-cM2} 
  For $f\in \cM^2_\omega(\cR)$, let
$f(0+i\cdot)=f_0$.  
Then $\cF f_0 \in L^2(\bR,e^{2e^\xi}d\xi)$ and
\begin{equation}\label{PW-cM2-equality}
\|f\|_{\cM^2_\omega(\cR)} = \|\cF f_0 \|_{L^2(\bR,e^{2e^{\xi}}d\xi)}\,.
\end{equation}

Conversely, if $\psi\in L^2(\bR,e^{2e^\xi}d\xi)$ and for $z\in\cR$
we set
\begin{equation}\label{PW-cM2-def}
f(z) =\frac{1}{\sqrt{2\pi}} \int_{-\infty}^{+\infty} \psi(\xi)
e^{z\xi}\, d\xi\, ,
\end{equation}
then $f\in \cM^2_\omega(\cR)$, equality\eqref{PW-cM2-equality} holds, and
$\psi=\cF f_0$.\ms
\end{theorem}

\begin{theorem}\label{repro-kernel}
The reproducing kernel for $\cM^2_\omega(\cR)$ is given by
$$
K(z,w)= \frac{1}{2\pi}\frac{\Gamma (z+\ov w)}{2^{z+\ov w}} \,. \ms
$$
\end{theorem}

We are in the position to describe the image of 
$A^2(\Delta)$ under the Mellin--Bergman transform $M_\Delta$.  Given 
any set $\Omega\subseteq\bC$ we denote by $\Hol(\Omega)$ the holomorphic
functions on $\Omega$. 
\begin{definition}{\rm
We define
\begin{equation}\label{H}
\cH = \big\{ g\in \Hol(\cR):\,  
\textstyle{\frac{\Gamma(1+z)}{2^z}} g(z) \in \cM^2_\omega(\cR) \big\}
\end{equation}
with norm
\begin{align*}
\|g\|_\cH^2 & = \Big\| \frac{\Gamma(1+z)}{2^z}  g
\Big\|_{\cM^2_\omega(\cR)}^2  
= 
\sum_{n=0}^{+\infty} \int_{-\infty}^{+\infty} |g({\tnt}+iy)|^2
\frac{|\Gamma({\tnt}+1+iy)|^2}{\Gamma(n+1)} \, dy \,. 
\end{align*}
}\ms
\end{definition}

\begin{theorem}\label{new}
The Mellin--Bergman transform
$$
M_\Delta: A^2(\Delta)\to \cH
$$ 
is a surjective isometry.
The space $\cH$ consists of holomorphic 
functions on $\cR$ that
are of exponential type at most $\pi/2$ and the polynomials are dense
in $\cH$.  Moreover, it is a Hilbert space with
reproducing kernel
$$
H(z,w) = \frac{1}{2\pi} \frac{\Gamma(z+\ov w)}{\Gamma(1+z)
\Gamma(1+\ov  w)} \, . 
$$
\end{theorem}
\ms

We stress the fact that we prove that
$$
M_\Delta\big(A^2(\Delta)\big) = \frac{2^z}{\Gamma(1+z)}
\cM_\omega^2(\cR)
=:\cH\, .
$$

It is interesting to notice that the space $\cH$ has already appeared
in the literature, in a different context \cite{KT2, KT1}\footnote{We
are  grateful to A. Aleman for pointing these references to us.}.  While on
one hand $\cH$ may be more natural being  the isometric image of $A^2(\Delta)$ 
through $M_\Delta$, the space $\cM^2_\omega(\cR)$ turns out to enjoy
more manageable properties, since the measure $\omega$ is translation
invariant in $\cR$ and there is no analogue of the Paley--Wiener type
Theorem \ref{PW-thm-cM2}  for $\cH$. 
We will collect some properties of $\cH$
and further remarks in Section \ref{Mellin}.  \ms

It is worth pointing  out that the measure $\omega$ was constructed in a
direct way and satisfies the uniqueness property in the next result.

\begin{theorem}\label{uniqueness-measure}
The measure $\omega$ is the unique positive translation invariant
Borel measure in $\ov\cR$
such that, for every $f,g\in A^2(\Delta)$ we have
$$
\la f,g\ra_{A^2(\Delta)}
= \big\la \textstyle{\frac{\Gamma(z+1)}{2^z}} M_\Delta(f), \,
  \textstyle{\frac{\Gamma(z+1)}{2^z}} M_\Delta(g) \big\ra_{L^2(\ov\cR,d\omega)}\,.
$$
\ms\end{theorem}

Clearly, the sets of uniqueness of
$\cH$ and $ \cM_\omega^2(\cR)$
coincide, and the same holds for the zero-sets.
We obtain that
 the zero-sets for $\cM^2_\omega(\cR)$ are
also zero-sets for the functions of exponential type $\pi/2$ and that
 zero-sets for the functions of exponential type $\tau<\pi/2$ are
 zero-sets  for $\cM^2_\omega(\cR)$, see Proposition
 \ref{zero-sets-containtment}. The reverse inclusions hold
 for the sets of uniqueness of the corresponding spaces.
\ms

In order to describe our next result we need to recall some classical
definitions.  
Given a sequence of points $\{z_j\}$ with $|z_j|\to+\infty$, 
 its exponent of convergence is 
$\rho_1= \inf \{\rho>0:\,  \sum_{j=1}^{+\infty} 1/|z_j|^\rho
<\infty\}$, while the counting function is 
$n(r)  = \# \{ z_j:\, |z_j|\le r\}$.  The upper and lower densities
$d^\pm=d_{\{z_j\}}^\pm$ are then defined as
$$
d^+ =  \limsup_{r\to+\infty} \frac{n(r)}{r^{\rho_1}}\,, \quad
d^- =  \liminf_{r\to+\infty} \frac{n(r)}{r^{\rho_1}}\,.
$$
\ms

In order to avoid  vanishing of infinite order at finite
points, we assume the functions to be regular in $\ov\cR$ and we
denote by $\Hol(\ov\cR)$ such space.
\begin{theorem}\label{zero-set-thm}
Let $\{z_j\}\subseteq\cR$, $1\le|z_j|\to+\infty$.  The following properties hold.
\begin{itemize}
\item[(i)] 
If $\{z_j\}$ has exponent of convergence $1$
 and  
upper density  $d^+<\frac12$, then
$\{z_j\}$ is a zero-set for $\cM^2_\omega(\cR)\cap\Hol(\ov\cR)$.\smallskip
\item[(ii)]
If $\{z_j\}$ is a zero-set for $\cM^2_\omega(\cR)\cap\Hol(\ov\cR)$, then
\begin{equation}\label{our-Carleman-cond}
\limsup_{R\to+\infty} \frac{1}{\log R} \sum_{|z_j|\le R} \Re\big(1/z_j\big) \le
\frac2\pi  \,.
\end{equation}
\end{itemize}
\end{theorem}

\begin{theorem}\label{MS-thm}
A sequence $\{z_j\}$ of points in $\cR$ such that $\Re z_j\ge \eps_0$,
for some $\eps_0>0$ and 
that violates condition \eqref{our-Carleman-cond}, is a set of
uniqueness for $\cM^2_\omega(\cR)$.

As a consequence, if $\{z_j\}$ is a sequence as above, the set of
powers $\{ \z^{z_j-1}\}$ is a complete set in $A^2(\Delta)$.  
\end{theorem}

We point out that a classical result of W. Fuchs' \cite{Fuchs} shows that 
there exist  sequences
 $\{z_j\}$ with exponent of convergence $1$ and lower density $d^->\frac12$ such
 that 
$\{z_j\}$ is a set of uniqueness for $\cM^2_\omega(\cR)\cap
\Hol(\ov\cR)$.  We will prove the above theorem and compare it with
Fuchs' result in Section \ref{MS-sec}.

\ms

For $1\le p<\infty$ we also consider the spaces
\begin{equation}\label{Mp-def}
\cM^p_\omega(\cR)
=\big\{ f\in \Hol(\cR): f\in H^p(S_b), \text{\ for all\ } b>0,\
\text{and\ }  f\in L^p(\ov\cR,d\omega) \big\}
\end{equation}
with norm
$$
\|f\|_{L^p(\ov\cR,d\omega)}^p =
\sum_{n=0}^{+\infty} \frac{2^n}{n!} \int_{-\infty}^{+\infty}
|f\big(\tnt +iy\big)|^p\, dy <\infty \,.
$$

We recall that $H^p(S_b) =\big\{ f\text{ holomorphic in } S_b:\,
  \sup_{0<x<b} \int_{-\infty}^{+\infty} |f(x+iy)|^p\, dy
  <\infty\big\}$ and that $f\in H^p(S_b)$ admits boundary values that
  are $p$-integrable, so that, as in the case $p=2$, integration over
  the imaginary is well defined.

Finally we prove 
\begin{theorem}\label{proj-thm}
The orthogonal projection operator 
$P: L^2(\ov\cR,d\omega) \to \cM^2$,
 is
unbounded as operator
$$
P: L^p(\ov\cR,d\omega)\cap L^2(\ov\cR,d\omega) \to \cM^p
$$
for every $p\neq2$.\ms
\end{theorem}

To the best of our knowledge, this is the first paper that deals with
a mixed Hardy--Bergman type condition that appears in the definition
of $\cM^2_\omega(\cR)$.  We find it remarkable that this
space appears naturally in the attempt of solving the \MS\ problem for the
Bergman space.  In \cite{Fuchs2}  Fuchs studied the \MS\ problem for
sets of exponential on the positive half-line.  
In \cite{Bologna} we study some generalization of
$\cM^2_\omega(\cR)$, obtained by  different choices of
the measure $\omega$.\ms

In order to relate ours to some previous work, we mention
that in \cite{Jacob-Partington-Pott1} 
B. Jacob, J. Partington and S. Pott
studied 
spaces of holomorphic function in $\cR$ whose norm is defined by the
condition 
$\sup_{\eps>0} \|f(\eps+\cdot)\|_{L^p(\ov\cR,d\mu)}<+\infty$, 
where $\mu$ is a translation invariant Borel measure on $\ov\cR$.
While, on one hand, this class of spaces contains as particular cases
the classical Hardy and Bergman spaces, $\cM^2_\omega(\cR)$ does not
fall in this class.  For, the finiteness 
of the above norm requires the function to be bounded in each
half-plane $\{ \Re z\ge \eps_0\}$, for $\eps_0>0$, 
while both $\cM^2_\omega(\cR)$ and $\cH$ contain 
functions of exponential growth.

We also mention that
in \cite{Sed} A. Sedletkskii studies the completeness of sets of
exponentials in weighted $L^p$ spaces on $(0,+\infty)$
in terms of zeros of functions the classical Bergman space on a
half-plane.  

\ms

\section{Basic properties of $\cM^2_\omega(\cR)$}\label{sec-cM2}
\ms

We begin recalling some well-known facts about Hardy spaces on a
strip.  For $0\le a<b<\infty$ we denote by $S_{(a,b)}$ the vertical strip
$\{z=x+iy: a<x<b\}$ and by $S_{[a,b]}$ its closure.  As before, we
simply write $S_b$ to denote the strip  $S_{(0,b)}$.  
The classical Hardy space
$H^2(S_{(a,b)})$ is 
$$
H^2(S_{(a,b)}) =\big\{ f\text{ holomorphic in } S_{(a,b)}:\,
  \sup_{a<x<b} \int_{-\infty}^{+\infty} |f(x+iy)|^2\, dy
  <\infty\big\}\,.
$$

\begin{thm}\label{PW-thm-strip}  {\rm {\bf (Paley--Wiener)}}

\noindent
{\rm (i)} Let  $F\in H^2(S_b)$.  Then $F_x:=F(x+i\cdot)$ admits limit
in $L^2(\bR)$ as $x\to0^+$ and $x\to b^-$, that we denote by $F_0$ and
$F_b$, respectively.  Moreover,  $e^{b\xi}\cF F_0 \in L^2(\bR)$
and for every $x\in[0,b]$
\begin{equation}\label{exp-form-PW}
 \cF F_x(\xi) = e^{x\xi} \cF F_0(\xi)\, .
\end{equation}

\noindent
{\rm (ii)}  $H^2(S_b)$ is a Hilbert space with the unique inner
product such that
$$
\|F\|_{ H^2(S_b)}^2 
= \|F_0\|_{L^2(\bR)}^2 + \|F_b\|_{L^2(\bR)}^2\,.
$$
{\rm (iii)} If $\psi\in L^2(\bR)$ and $e^{b\xi} \psi \in L^2(\bR)$,
then 
$$
F(z) =\frac{1}{\sqrt{2\pi}} \int_{-\infty}^{+\infty}
\psi(\xi)e^{z\xi}\, d\xi
$$
is in $H^2(S_b)$ and $\cF F_0=\psi$.
\end{thm}

The proofs of these facts can be found in \cite{PW}, Theorems
I-VII. \ms

Next we prove Theorem \ref{cM2-Hilbert}.  From now on, for simplicity
of notation, we write $\cM^2$ in place of $\cM^2_\omega(\cR)$.

\proof[Proof of Theorem \ref{cM2-Hilbert}]
 This is elementary and we include the details for
completeness.
 
Let $\{f_m\}$
be  sequence in  $\cM^2$, Cauchy in the $L^2(\ov\cR,d\omega)$-norm. 
Then $\{f_m(\frac{k}{2}+i\cdot)\}$ is a Cauchy sequence in $L^2(\bR)$
for every $k$, since
$$
\| f_{m}-f_{m'}\|_{L^2(\ov\cR,d\omega)}^2
= \sum_{k=0}^{+\infty} \frac{2^k}{k!} \int_{-\infty}^{+\infty}
|f_m\big(\tkt +iy\big) -f_{m'}\big(\tkt +iy\big)|^2\, dy\\
<\eps
$$
implies that for any fixed $n$,
$$
\|f_m\big(\tkt +i\cdot\big) -f_{m'}\big(\tkt +i\cdot\big)\|_{L^2(\bR)}^2
<
\eps' 
$$
for $m,m'\ge N$, and $0\le k\le n$. 
By (ii) in the previous theorem, $\{f_m\}$ is 
a Cauchy sequence in $H^2(S_{\frac{n}{2}})$, so that  
$\{f_m\}$ converges to $f \in H^2(S_{\frac{n}{2}})$, for any $n$.  By analytic
continuation $f\in H^2(S_b)$ for all $b>0$.  

Now it is clear that $f_m\to f$ in the $L^2(\ov\cR,d\omega)$-norm, that
is, $\cM^2$ is closed.  
The fact that it is a Hilbert space with
reproducing kernel follows at once, since point evaluations are
bounded in $H^2(S_b)$ and that if $f\in \cM^2$, its
$H^2(S_b)$-norms are controlled by a constant (depending on $b$)
times the $L^2(\ov\cR,d\omega)$-norm.
The conclusion about the inner product is now clear.  
\qed
\ms

We take a look at the elementary inclusions between the basic function
spaces. 
\begin{prop}\label{prop-1.2}
\ 
\begin{itemize}
\item[(i)] Let $0<\delta \le1$ and $\eps_0>0$. Then 
 $\Gamma(\eps_0+\delta z)\in \cM^2$ if and only if $\delta <1$.\smallskip
\item[(ii)] Let $G(z) = (z+1)^{-1} \exp\big\{ ie^{2\pi iz}\big\}$. Then 
$G\in L^2(\ov\cR, d\omega)$ is  holomorphic in $\ov\cR$,  but $G\not\in
  H^2(S_b)$ if $b>\frac12$. On the other hand,
  $\Gamma(1+z)\in H^2(S_b)$ for $b>0$, but
  $\Gamma(1+z)\not\in \cM^2$. 
\item[(iii)] Let $h$ be a function regular in $\ov\cR$ and of
  exponential type $\tau<\pi/2$ and let $2\tau/\pi<\delta <1$.  Then 
 $F(z)=h(z)\Gamma(1+\delta z)\in \cM^2$.
\end{itemize}
\end{prop}

\proof
It is well known that  
for $c>0$,
\begin{align*}
\|\Gamma(c+i\cdot)\|_{L^2(\bR)}^2 
& = \int_0^{+\infty} e^{-2x} x^{2c}\, \frac{dx}{x} = 2^{-2c} \Gamma(2c)\,,
\end{align*}
see Lemma 2.3 in \cite{BuJa} e.g. (and also
the discussion in Section \ref{Mellin}).  Therefore, for every $b>0$,
$$
\sup_{0<x<b} \| \Gamma(\eps_0+\delta x +i\delta\cdot)\|_{L^2(\bR)}^2 
= C_b<+\infty\, .
$$

Thus, in order to prove $\Gamma(\eps_0+\delta z)\in
\cM^2$ it suffices to show that $\Gamma(\eps_0+\delta z)\in
L^2(\ov\cR,d\omega)$.  

For some $Q>0$ large enough we have
\begin{align*}
\|\Gamma(\eps_0+\delta \cdot)\|_{\cM^2}^2 
& = \sum_{n=0}^{+\infty} \frac{2^n}{n!} 
\|\Gamma(\eps_0+\delta \tnt +i\delta \cdot)\|_{L^2(\bR)}^2 \\
& \le\frac1\delta  \sum_{n=0}^{+\infty} \frac{2^{(1-\delta )n}}{n!}
\Gamma(2\eps_0+\delta n) \\
& \le C \sum_{n=0}^{+\infty} \frac{2^{(1-\delta )n}}{n!} n^{\eps_0}
  \delta ^n e^{\delta n\log n} \\
& \le C \sum_{n=0}^{+\infty} Q^n
 e^{(\delta -1)n\log n} \\
& <\infty\,.
\end{align*}
Since it is easy to see that the norm
$\|\Gamma(\eps_0+\delta \cdot)\|_{\cM^2}$ is infinite
 if $\delta =1$, this proves (i).

In order to prove (ii), notice that 
$$
\big| \exp\big\{ ie^{2\pi iz}\big\}\big| = \exp\big\{ e^{-2\pi
  y}\cos(2\pi x) \big\}\,,
$$
is bounded for $x=\tnt$, $n=0,1,2\dots$.  This implies
that $G(z) = (z+1)^{-1} \exp\big\{ ie^{2\pi iz}\big\} \in
L^2(\ov\cR,d\omega)$.  However, $G\not\in H^2\big(S_b\big)$ if
$\cos(2\pi x)>0$ for some $x<b$. \ms

Finally, recall the asymptotic of the Gamma function (see \cite{Lebedev}
  e.g.), valid 
for $|\arg z|\le \pi-\delta$,
\begin{equation}\label{asy-basic}
\Gamma(z) =\sqrt{2\pi}\, e^{(z-\frac12)\log z-z}\Big[ 1+\cO(1/|z|)\Big]
\,.
\end{equation}
Then $\Gamma(\eps_0+\delta z)$ is regular in $\ov\cR$ and for
$|z|\ge1$, $z=x+iy$, 
\begin{align}
|\Gamma(\eps_0+\delta z)|
& \le C  |z|^{-\frac12}
\exp \big\{ \Re [(\eps_0+\delta z)\big(\log(\eps_0+\delta z) -1\big)]\big\} 
\notag \\
& \le C |z|^{\eps_0-\frac12}\exp \big\{ \delta x\log|z| -\delta |y|\arctan
(|y|/x)\big\}\,. \label{Gamma-fnc-est}
\end{align}

  Let $2\tau/\pi<\delta <1$ and take
$\eps_0=1$ (for simplicity).  When $0\le x\le b$ we have
$\tau < \tau'<\delta \arctan (|y|/x)$ for $|y|\ge N$ sufficiently large.
Then \eqref{Gamma-fnc-est} gives
$$
|h(z)|\,| \Gamma(1+\delta z)| \le e^{\tau|z|} e^{-\tau'|y|} e^{\delta x\log
  |z|}\le C
$$
uniformly in the strip $S_b$.  In order to bound the
$\cM^2$-norm we observe that $|h(z)|\le C e^{\tau(x+|y|)}$ and
using \eqref{Gamma-fnc-est} we estimate
\begin{align*}
& \int_{-\infty}^{+\infty} e^{2\tau|y|}  |\Gamma(1+\delta \tnt +i\delta y)|^2\,
dy \\ 
& \qquad\le \int_{-\infty}^{+\infty} e^{2\tau|y|}  (n+|y|) \exp\big\{
 \delta n \log\big(\tnt+|y|\big) -2\delta |y|\arctan \big(2|y|/n\big)\big\}  \, dy \\ 
& \qquad = \bigg( \int_{|y|\le \alpha n} + \int_{|y|>\alpha n} \bigg)   (n+|y|) \exp\big\{2\tau|y|+
 \delta n \log\big(\tnt+|y|\big) -2\delta |y|\arctan \big(2|y|/n\big)\big\}  \, dy
 \\
& \qquad = I +I\!I  \,,
\end{align*}
with $\alpha >\frac12$ to be fixed.  Then it follows at once that 
\begin{equation}\label{est-I}
I \le C n^2 (1+\alpha)^{\delta n}e^{2\tau\alpha n} e^{\delta n\log n}\le C p^n e^{\delta n\log n}\,,
\end{equation}
for some $p>0$.
On the other hand, since $\tau/\delta<\pi/2$ we can select $\alpha$
large enough so that
$\delta\arctan(2\alpha)>\tau$.  Then we have
\begin{align}
I\!I 
& \le C \int_{|y|>\alpha n}  |y| \exp\big\{ 2\tau|y|+
 \delta n \log\big(\tnt+|y|\big) -2\delta |y|\arctan \big(2|y|/n\big)\big\}  \, dy
\notag \\
& \le C 2^{\delta n} \int_{|y|>\alpha n}  |y|^{1+\delta n}  
\exp\big\{ 2|y|\big(\tau -\delta \arctan \big(2|y|/n\big)\big)\big\}
\, dy \notag \\
& \le C 2^{\delta n} \int_{|y|>\alpha n}  |y|^{1+\delta n}  
\exp\big\{ -2\big(\delta \arctan (2\alpha)-\tau\big)|y|\big\}
\, dy  \notag \\
& \le C 2^{\delta n} \frac{\Gamma(2+\delta n)}{ 
\big[ 2\big(\delta \arctan (2\alpha)-\tau\big)\big]^{2+\delta n}}
\notag \\
& \le C q^n e^{\delta n\log n}
\,,\label{est-II}
\end{align}
for some $q>0$ and where we have applied the simple estimate, valid
for $A,B,y_0>0$,
$$
\int_{y_0}^{+\infty} e^{-Ay}y^B\, dy 
= \frac{1}{A^{B+1}} \int_{Ay_0}^{+\infty} e^{-t} t^B\, dt \le
\frac{\Gamma(B+1)}{A^{B+1}}\,.
$$ 

Putting  \eqref{est-I} and \eqref{est-II} together we then obtain
 \begin{align*}
\|h\, \Gamma(1+\delta \cdot)\|_{\cM^2}^2 
& \le C \sum_{n=0}^{+\infty} \frac{2^n}{n!} 
(p^n+ q^n) e^{\delta n\log n}\\
& \le C \sum_{n=0}^{+\infty} Q^n e^{(\delta-1) n\log n}\\
& < \infty\,.
\end{align*}
This proves (iii).  
\qed
\ms

\proof[Proof of Theorem \ref{PW-thm-cM2}]
Since $f\in H^2(S_{\frac{n}{2}})$ for every $n$, using \eqref{exp-form-PW}
we have,
\begin{align}
\|f\|_{\cM^2}^2 
& = \sum_{n=0}^{+\infty} \frac{2^n}{n!} \| 
f(\tnt+i\cdot)\|_{L^2(\bR)}^2 \notag \\
& = \sum_{n=0}^{+\infty} \frac{2^n}{n!} \big\|
\cF\big(f(\tnt+i\cdot)\big)\big\|_{L^2(\bR)}^2 \notag \\ 
& = \sum_{n=0}^{+\infty} \frac{2^n }{n!} \int_{-\infty}^{+\infty}
e^{n\xi}|\cF f_0(\xi)|^2 \, d\xi \notag \\
& = \int_{-\infty}^{+\infty}
|\cF f_0(\xi)|^2 e^{2e^{\xi}}\, d\xi  \label{one-more-PW}
\,.
\end{align}

Conversely, let $\psi\in L^2(\bR,e^{2e^\xi}d\xi)$ and $f$ be defined by
\eqref{PW-cM2-def}.  Notice that the integral converges absolutely for
$z\in\cR$ since
$$
\int_{-\infty}^{+\infty} |\psi(\xi) e^{z\xi}|\, d\xi 
\le \|\psi\|_{L^2(\bR,e^{2e^\xi}d\xi)}  \bigg(\int_{-\infty}^{+\infty}
e^{2x\xi} e^{-2e^\xi}\, d\xi \bigg)^{1/2}<\infty\,.
$$
Then $f$ is well defined and holomorphic in $\cR$.  Notice also that
$e^{\frac{n}{2}\xi}\psi\in L^2(\bR)$ for $n=0,1,2,\dots$
so that by (iii) in Theorem \ref{PW-thm-strip}
$f\in H^2(S_{\frac{n}{2}})$ for every $n$,  and
$\cF f_0=\psi$.  Now
the same argument as in \eqref{one-more-PW} gives \eqref{PW-cM2-equality}.
\qed
\ms

As consequence of the Paley--Wiener-type theorem we determine the
reproducing kernel of $\cM^2$.

\proof[Proof of Theorem \ref{repro-kernel}]
Let $K_z\in \cM^2$ be such that $\la f,\,
K_z\ra_{\cM^2}=f(z)$ for every $z\in\cR$ and every $f\in
\cM^2$.   
Using \eqref{exp-form-PW} again we have,
\begin{align*}
f(z) 
& = \sum_{n=0}^{+\infty} \frac{2^n}{n!} \big\la f(\tnt+i\cdot),\,
K_z(\tnt+i\cdot)\big\ra_{L^2} \\
& = \sum_{n=0}^{+\infty} \frac{2^n}{n!} \big\la \cF\big(f(\tnt+i\cdot)\big),\,
\cF\big(K_z(\tnt+i\cdot)\big)\big\ra_{L^2} \\
& = \sum_{n=0}^{+\infty} \frac{2^n }{n!} \big\la e^{\frac{n}{2}\xi}\cF f_0,\,
e^{\frac{n}{2}\xi} \cF K_{z,0} \big\ra_{L^2} \\
& = \int_{-\infty}^{+\infty} \cF f_0(\xi)
\ov{\cF K_{z,0}(\xi)}\, e^{2e^\xi}\, d\xi\,,
\end{align*}
where switching the integral with the sum is justisfied by Theorem
\ref{PW-thm-cM2}. 

On the other hand, 
$$
f(z) =\frac{1}{\sqrt{2\pi}} \int_{-\infty}^{+\infty} e^{z\xi}\cF
f_0(\xi)\, d\xi\,,
$$
so that 
$$
\cF K_{z,0} (\xi) = \frac{1}{\sqrt{2\pi}} e^{-2e^{\xi}} e^{\ov z\xi} \,,
$$
and
\begin{align*}
K_z(w) 
&= \frac{1}{2\pi} \int_\bR e^{w\xi}  e^{-2e^{\xi}} e^{\ov z\xi}\, d\xi\\
& = \frac{1}{2\pi} \int_\bR e^{(w+\ov z) \xi}  e^{-2e^{\xi}}\, d\xi \\
&= \frac{1}{2\pi 2^{w+\ov z}}  \int_0^{+\infty} t^{w+\ov z-1}
e^{-t} \, dt \\
& = \frac{1}{2\pi} \frac{\Gamma(w+\ov z)}{2^{w+\ov z}} \,. \qed
\end{align*}
\ms

\section{Properties of  Mellin--Bergman transforms}\label{1}
\ms

Next we study the Mellin--Bergman transform $M_\Delta $, as defined in 
\eqref{M-Bt}, when acting on 
$A^2(\Delta)$. We set $\|f\|_{A^2(\Delta)}^2 = \frac{1}{\pi}\iint_\Delta |f(\z)|^2\, dA(\z)$.
Observe that $M_\Delta f(z) =\la f,\, \z^{\ov z-1}\ra_{A^2(\Delta)}$,
so that $M_\Delta$ is linear and $M_\Delta f$ is
holomorphic in $\cR$.

We need the explicit expression of the inner product of the powers 
$\z^\alpha$ and $\z^\beta$ in $A^2(\Delta)$.  The next result appears
in \cite{KPS2}.
\begin{lem}\label{Beta-fnc} 
Let $\Re\alpha,\Re\beta>-1$.  Then
$$
\frac{1}{\pi} \iint_\Delta \z^\alpha\overline{\z^\beta}\, dA(\z)
=  \frac{\Gamma(\alpha+\overline\beta+2)}{
\Gamma(\alpha+2)\Gamma(\overline\beta+2)}\, .
$$
In particular, $\z^{z-1}\in A^2(\Delta)$ if and only if
$\Re z>0$ and in this case
$$
\|\z^{z-1}\|^2_{A^2(\Delta)} = \|\z^{\ov z-1}\|^2_{A^2(\Delta)}= 
\frac{\Gamma(2\Re z)}{|\Gamma(z+1)|^2}\, .
$$
\end{lem}

The lemma provides us with a few explicit examples.
\begin{exam}\label{imag-polys}{\rm 
By taking $\alpha=0$ and
$\beta=\ov z-1$, we obtain that
$M_\Delta(1) = 1$, and similarly, for $k=1,2,\dots$
\begin{align*}
M_\Delta (\z^k) (z)
& = \ \frac{\Gamma(z+k+1)}{k! \,\Gamma(z+1)}\\
& = \frac{1}{k!}\, (z+1)\cdots(z+k)\, ;
\end{align*}
hence a polynomial of degree $k$ in $z$.}
\end{exam}

We now obtain the following  property on the growth of functions in
$M_\Delta\big(A^2(\Delta)\big)$. 
\begin{prop}\label{Mf-exp-type}
For $f\in A^2(\Delta)$, 
$M_\Delta f$ is holomorphic of exponential type at most $\frac\pi2$
in $\cR$. 
\end{prop}

\proof 
It is clear that $M_\Delta f$ is holomorphic in $\cR$ since
\begin{align*}
\big| \p_z \big( f(\z) \ov \z^{z-1}\big) \big|
&\le c | f(\z)|\, \big|  \log|\z|\big| \big| \ov \z^{z-1}\big| \\
& \le  c |f(\z)| \, |\z|^{\Re z  -\eps-1} e^{\Im z \pi/2}\, ,
\end{align*}
which is absolutely integrable if $\Re z>\eps$. \ms

Next, 
\begin{equation}
|M_\Delta f(z)| 
 \le  \|f\|_{A^2(\Delta)} \|\z^{\ov z -1}\|_{A^2(\Delta)}   \le C \frac{\Gamma(2\Re
  z)^{1/2}}{|\Gamma(z+1)|}\,.  \label{tag}
\end{equation}
A straightforward application of 
the asymptotics of the Gamma function \eqref{asy-basic} shows that 
$M_\Delta f$ is of exponential type at most $\pi/2$.  
We leave the details to the reader.

We also observe that there exist functions $f$ such that $M_\Delta f$
is of type $\pi/2 - \eps$, for every $\eps>0$. For, for $\z_0\in
\Delta$, let $B_{\z_0}$ be the reproducing kernel for $A^2(\Delta)$ at
$z_0$.  Then,
\begin{align*}
|M_\Delta(B_{\z_0})(z)|
& =  |\ov\z_0^{z-1}| \ge e^{y\arg\z_0}\, . 
\end{align*}
The conclusion follows by taking $\z_0$ with $\arg\z_0\ge \pi/2-\eps$
and $z\in\cR$, $z=\delta+iy$.
\qed\ms

We now turn to proving Theorem \ref{new},  therefore
establishing  the
properties the Mellin--Bergman transforms of
$A^2(\Delta)$-functions.  

We break  Theorem \ref{new}  into two results, first of which is the
following theorem.
\ms

 As usual, we denote by $A^2(\cR)$ the (unweigthed)
Bergman space, then the Paley--Wiener theorem for $A^2(\cR)$
shows that the
Fourier transform is a surjective isometry between
$A^2(\cR)$ and $L^2\big((-\infty,0), d\xi/|\xi|\big)$. 

\begin{thm}\label{PW-thm} 
Let  $F\in A^2(\cR)$. Then there exists $\psi\in L^2 \big((-\infty,0),
d\xi/|\xi|\big)$
such that for $z\in\cR$
\begin{equation}\label{PW-ide}
F(z) =\frac{1}{\sqrt{2\pi}} \int_{-\infty}^0 e^{z \xi} \psi(\xi)\, d\xi
\end{equation}
and
\begin{equation}\label{PW-iso}
\|F\|^2_{A^2(\cR)}  = \textstyle{\frac12}\|\psi\|^2_{L^2 ((-\infty,0),\, d\xi/|\xi|)} \,.
\end{equation}
Conversely, if $\psi\in L^2\big((-\infty,0), d\xi/|\xi|\big)$ and $F$ is
defined by \eqref{PW-ide} then $F\in A^2(\cR)$ and \eqref{PW-iso} holds.
\end{thm}

A proof of this result can be found in \cite{Be-al}  or \cite{Duren-et-al}, e.g. 
This result implies that if $F\in A^2(\cR)$, setting
\begin{equation}
F_0=\lim_{x\to 0^+} \cF^{-1}\big(e^{x(\cdot)}\psi\big)
\end{equation}
in the sense of tempered distributions, we obtain that $F$ admits
boundary values $F_0=F(0+i\cdot)$ such that $\cF F_0= \psi$.

The next result can be
obtained from the well-known property of the Gamma function
for $\Re \lambda >-1$,
\begin {equation}
\int_0^{+\infty} t^{\lambda}e^{-t}e^{ixt}\,dt=\frac{\Gamma(\lambda +1)}
{(1-ix)^{\lambda +1}} \,.
\end{equation}

\begin{lem}\label{a-Fourier-transform}
Let $\Re\lambda>0$ and  $h(w)
=(1+w)^{-\lambda-1}$,  then $h\in A^2(\cR)$ and
$$
\cF h_0 (\xi)= \sqrt{2\pi}
\frac{|\xi|^\lambda}{\Gamma(1+\lambda)}
e^{\xi}\chi_{\{\xi<0\}} \,. \ms
$$
\end{lem}

We recall that the (re-normalized) classical Mellin transform is 
\begin {equation}\label{Mellin-def}
M\vp (z) = \frac{1}{\sqrt{2\pi}} \int_0^{+\infty} \vp(t) t^{z-1}\, dt\,,
\end{equation}
where $\vp$ is a function defined on $(0,+\infty)$. 

Next, we reduce the problem to characterize the space
$M\big(L^2 \big(  (0,+\infty),\, \frac{e^{2\xi}}{\xi}\,
d\xi\big)$.  

\begin{thm}\label{reduction-thm}
There exists a surjective isometry
$$
T: A^2(\Delta)\to L^2 \Big( (0,+\infty),\, \frac{e^{2\xi}}{\xi}\,
d\xi\Big)
$$
such that for $g\in A^2(\Delta)$ and $z\in\cR$
$$
M_\Delta g(z) = -\sqrt{\pi} \frac{2^z}{\Gamma(z+1)} M(Tg)(z)\,.
$$
\end{thm}

\proof
 For $w \in\cR$
we let  $\phi(w) =2(w+1)^{-1}$.   Then $\phi:\cR\to\Delta$ is
a biholomorhic mapping and $f\mapsto\frac1{\sqrt \pi} \phi' (f\circ\phi)=:\tilde f$ is
a surjective isometry of $A^2(\Delta)$ onto $A^2(\cR)$.  Notice that
$$
\big( \zeta^{\ov z-1}\big)\widetilde{\ } (w) = - \frac{2^{\ov z}}{\sqrt \pi}
\frac{1}{ (w+1)^{\ov z+1}}\,.
$$

Then, if $f\in A^2(\Delta)$, using Lemma \ref{a-Fourier-transform}
(with $z=\ov\lambda\in\cR$) 
we have
\begin{align*}
M_\Delta f(z) 
& = \frac 1 \pi \iint_\Delta f(\z) \ov \z^{\ov z-1}\, dA(\z) = \la f,\,  \z^{\ov z-1} \ra_{ A^2(\Delta)} \\
& = -\frac{2^{z}}{\sqrt \pi} \big\la \tilde f,\,  (w+1)^{-\ov z-1}
\big\ra_{A^2(\cR)} \\
& = -\sqrt{\frac12} \frac{2^z}{\Gamma(1+z)} \big\la \cF\tilde f_0,\,
 (- \xi)^{\ov z} e^{\xi} \big\ra_{L^2((-\infty,0),\, d\xi/|\xi|)} \\
& = -\sqrt{\frac 12} \frac{2^z}{\Gamma(1+z)}
\int_{-\infty}^0\cF \tilde f_0(\xi) e^{\xi} (-\xi)^{z-1}\, d\xi \\
& =  -\sqrt{\pi} \frac{2^z}{\Gamma(1+z)} M\big(
e^{-t} \cF\tilde f _0(-t) \big) (z)\,;
\end{align*}
$M$ being the Mellin transform.  Setting for $t>0$, $Tf(t )=e^{-t} \cF\tilde f_0 (-t)$
the desired conclusion follows at once.
\qed\ms

\ms

\section{Mapping properties of the Mellin transform}\label{Mellin}
\ms

In order to complete the proof of Theorem \ref{new} we need the
following result.

\begin{thm}\label{mapping-M}
The mapping
$$
M: L^2 \Big( (0,+\infty),\, \frac{e^{2\xi}}{\xi}\, d\xi\Big) \to
\cM^2
$$
is a surjective isometry.
\end{thm}

The mapping properties of $M$ as operator between
function spaces have been studied in \cite{BuJa}.

We begin the proof of theorem  with the following
\begin{lem}\label{Mellin-unitary}
The mapping
$$
M: L^2 \Big( (0,+\infty),\, \frac{e^{2\xi}}{\xi}\, d\xi\Big) \to 
\cM^2
$$
is a partial isometry.
\end{lem}

\proof
We first show that if $\vp\in L^2 \big( (0,+\infty),\,
\frac{e^{2\xi}}{\xi}\, d\xi\big)$ then $M\vp$ is well defined and
holomorphic in $\cR$.  Writing $z=x+iy$ we have 
\begin{align*}
|M\vp(z)|
& \le \frac{1}{\sqrt{2\pi}} \int_0^{+\infty} |\vp(\xi)| \xi^{x-1}\, d\xi\\
& \le \frac{1}{\sqrt{2\pi}} \| \vp \|_{L^2 ( (0,+\infty),\frac{e^{2\xi}}{\xi}\, d\xi)} \Big(
\int_0^{+\infty} e^{-2\xi}  \xi^{2x-1}\, d\xi \Big)^{1/2} \\
& = \frac{1}{\sqrt{2\pi}} 
\frac{\Gamma(2x)^{1/2}}{2^x} \| \vp \|_{L^2 ( (0,+\infty),\frac{e^{2\xi}}{\xi}\, d\xi)}
\,. 
\end{align*}
Hence, the integral defining $M\vp$ converges absolutely for
every $z\in \cR$, and a similar argument shows that $M\vp$ is also
holomorphic in $\cR$.  
Notice that  
\begin{align*}
M\vp(z) 
& = \frac{1}{\sqrt{2\pi}} \int_{-\infty}^{+\infty} \big(\vp\circ\exp)(s) e^{zs}\, ds\\
& = \cF^{-1} \big( (\vp\circ\exp) e^{x(\cdot)}\big)
(y)\,.
\end{align*}
Moreover, 
\begin{align*}
\big\| M\vp(x+i\cdot) \big\|_{L^2(\bR)}^2 
& = \| \vp \|_{L^2( (0,+\infty), \xi^{2x-1}d\xi)}^2 \\
& \le C_x  \| \vp \|_{L^2( (0,+\infty), \frac{e^{2\xi}}{\xi} d\xi)}^2
\, ,
\end{align*}
 uniformly in $x\in (0,b]$.
Hence,  $\vp\in L^2 \big( (0,+\infty),\, \frac{e^{2\xi}}{\xi}\,
d\xi\big)$ implies that   $M\vp \in H^2(S_b)$ for every $b>0$. \ms

Finally, let $\vp,\psi\in  L^2 \big( (0,+\infty),\,
\frac{e^{2\xi}}{\xi}\, d\xi\big)$ and first assume that they have compact
support.  Then,  we have
\begin{align*}
\la M\vp,\, M\psi\ra_{\cM^2}
& = \frac{1}{2\pi} \sum_{n=0}^{+\infty} \frac{2^n}{n!} \int_{-\infty}^{+\infty} 
M\vp (\tnt +iy) \ov{M\psi (\tnt +iy) }\, dy \\
& = \sum_{n=0}^{+\infty} \frac{2^n}{n!} \int_{-\infty}^{+\infty} 
\cF^{-1} \Big( \big(\vp\circ\exp) e^{\frac{n}{2}(\cdot)}\Big) (y)
\ov{ \cF^{-1} \Big( \big(\psi\circ\exp) e^{\frac{n}{2}(\cdot)}\Big) (y)} \, dy \\
& = \sum_{n=0}^{+\infty} \frac{2^n}{n!} \int_{-\infty}^{+\infty} 
\vp \big(e^y\big) e^{\frac{n}{2}y} 
\ov{\psi \big(e^y\big) } e^{\frac{n}{2}y}  \, dy \\
& = \int_{-\infty}^{+\infty} 
\vp \big(e^y\big) \ov{\psi \big(e^y\big) } e^{2e^y}\, dy \\
& = \int_0^{+\infty} \vp (t)
\ov{\psi (t) } \, e^{2t} \frac{dt}{t}\,. 
\end{align*}
Therefore, 
$ M: L^2 \big( (0,+\infty),\, \frac{e^{2\xi}}{\xi}\, d\xi\big) \to
\cM^2$ is a partial isometry, i.e.
$$
 \la M\vp,\, M\psi\ra_{\cM^2}
= \la \vp,\,\psi\ra_{ L^2( (0,+\infty),\,\frac{e^{2\xi}}{\xi}\,
    d\xi)}
$$
for all $\vp,\psi\in L^2\big( (0,+\infty),\,\frac{e^{2\xi}}{\xi}\,
    d\xi\big)$.  
\qed
\ms

To order to prove Theorem \ref{mapping-M} we need some density results
in $\cM^2$, that are consequences of the Paley--Wiener-type Theorem \ref{PW-thm-cM2}.

For $1\le p<\infty$ we 
 denote by $\cM^p_{(\eps)}$ the subspace of 
$\cM^p$ of functions that are
 holomorphic for $\Re z>-\eps$ and 
that are in $H^p(S_{(-\eps,b)})$ for every $b>0$.

\begin{lem}\label{density-1}
Let $\eps'>0$ and  $\psi \in L^2\big(\bR,
e^{2\eps'\xi}e^{2e^\xi}d\xi\big)\cap  L^2\big(\bR,
e^{2e^\xi}d\xi\big)$. Then
$\psi \in L^2\big(\bR,
e^{2\eps\xi}e^{2e^\xi}d\xi\big)$ for $0<\eps\le  \eps'$ and let
$f$ be defined by \eqref{PW-cM2-def}.  If $f_\eps(z)=f(z+\eps)$, then
$f_\eps\in \cM^2_{(\eps)} $ and 
$f_\eps \to f$ in $\cM^2$ as $\eps\to 0^+$.  Hence, $\bigcap_{\eps>0}\cM^2_{(\eps)} $ is
dense in $\cM^2$.
\end{lem}
\proof
It is clear that $\psi \in L^2\big(\bR,
e^{2\eps\xi}e^{2e^\xi}d\xi\big) \cap  L^2\big(\bR,
e^{2e^\xi}d\xi\big)$ for $0<\eps< \eps'$.
Notice that 
$$
f_\eps(z) = \frac{1}{\sqrt{2\pi}} \int_{-\infty}^{+\infty} \psi(\xi)
e^{(z+\eps)\xi}\, d\xi\,,
$$
it is holomorphic  in $\Re z>-\eps$, 
and since $e^{\eps(\cdot)}\psi \in L^2\big(\bR, e^{2e^\xi}d\xi\big)$,
  $f_\eps \in \cM^2$.  
Moreover, since $f\in H^2(S_b)$ for every $b>0$, $f_\eps \in H^2(S_{(-\eps,b)})$ for
every $b>0$.  Hence, $f_\eps \in \cM^2_{(\eps)}$. 

Next, since
$$
\cF\big( f(x+\eps+i\cdot)\big)(\xi) = \psi(\xi)
e^{(x+\eps)\xi}\,,
$$
we have
\begin{align*}
\|f-f_\eps\|_{\cM^2}^2
 & = \sum_{n=0}^{+\infty} \frac{2^n}{n!} \big\|
\cF\big(f(\tnt+i\cdot)\big) - \cF\big( f(\tnt+\eps+i\cdot)\big)
\big\|_{L^2(\bR)}^2 \notag \\ 
& = \sum_{n=0}^{+\infty} \frac{2^n}{n!} 
\int_{-\infty}^{+\infty}
e^{n\xi}|(1-e^{\eps\xi})\psi(\xi)|^2 \, d\xi\\
& = \int_{-\infty}^{+\infty}
|(1-e^{\eps\xi})\psi(\xi)|^2 e^{2e^{\xi}}\, d\xi \to 0\,, 
\end{align*}
as $\eps\to0$.  Since $L^2\big(\bR, e^{2\eps\xi}e^{2e^\xi}d\xi\big)\cap  L^2\big(\bR,
e^{2e^\xi}d\xi\big)$ is dense in 
$L^2\big(\bR, e^{2e^\xi}d\xi\big)$ the 
 conclusion follows.
\qed
\ms

\begin{prop}\label{density-2}
The subspace
$\bigcap_{\eps>0}\cM^2_{(\eps)}\cap \cM^1_{(\eps)}$ is dense in
$\cM^2$. 
\end{prop}
\proof
Let $\psi\in C^\infty_0(\bR)$ and $f$ be defined  by
\eqref{PW-cM2-def}.  By the previous result $f\in
\cM^2_{(\eps)} $.  Now, denoting by $\|\psi\|_{W^s(\bR)}$ the standard
Sobolev norm, 
\begin{align*}
\|f\|_{L^1(\cR,d\omega)}
& \le \|(1+y^2)^{-1/2}\|_{L^2(\ov\cR,d\omega)}
\|(1+y^2)^{1/2}f\|_{L^2(\ov\cR,d\omega)} \\
& \le C \bigg( \sum_{n=0}^{+\infty} \frac{2^n}{n!} \int_{-\infty}^{+\infty} (1+y^2)
\big|\cF^{-1}\big(e^{\frac{n}{2}(\cdot)}\psi\big)(y)\big|^2\,
dy\bigg)^{1/2} \\
& = C \bigg( \sum_{n=0}^{+\infty} \frac{2^n}{n!}
\big\|e^{\frac{n}{2}(\cdot)}\psi\big\|_{W^2(\bR)}^2 \bigg)^{1/2}\\
& \le C \bigg( \sum_{n=0}^{+\infty} \frac{2^n}{n!}
\int_{-\infty}^{+\infty} e^{n\xi} \big( n^4|\psi(\xi)|^2
+n^2|\psi'(\xi)|^2 + |\psi''(\xi)|^2 \big) d\xi\bigg)^{1/2} \\
&  \le C \bigg( 
\int_{-\infty}^{+\infty}  \big(|\psi(\xi)|^2
+n^2|\psi'(\xi)|^2 + |\psi''(\xi)|^2 \big) P(d\xi/\xi) 
\big( e^{2e^\xi}\big) d\xi\bigg)^{1/2}
\end{align*}
where $P$ is the polynomial $t^4+t^2+1$.  The right hand side above
is finite since $\psi\in C^\infty_0$, so that $f\in L^1(\ov\cR,d\omega)$.  
Arguing as before, we see that for every $x\ge-\eps$
$$
\|f(x+i\cdot)\|_{L^1(\bR)} \le C \| e^{x(\cdot)} \psi \|_{W^2(\bR)}^2
\, ,
$$
hence $f\in H^1(S_{(-\eps,b)})$ for every $b>0$, that is, $f\in
\cM^2_{(\eps)} \cap \cM^1_{(\eps)}$. 
Since 
$C^\infty_0$ is dense in $L^2\big(\bR,e^{2e^\xi} d\xi\big)$ the
conclusion follows.
\qed
\ms

\proof[Proof of Theorem \ref{mapping-M}] We only need to
show that $M$ is onto.
\ms

It is a well-known fact  that if  
$g\in L^2\big(\{c\}+i\bR\big) \cap
L^1 \big(\{c\}+i\bR\big)$ for all $c\in(a,b)$ and
$g(x+iy)\to 0$ as $|y|\to+\infty$ uniformly in $x\in(a,b)$, then for $\xi>0$
\begin{equation}\label{M-inv-c}
M_c^{-1} g(\xi) =\frac{1}{\sqrt{2\pi}} 
\int_{-\infty}^{+\infty} g(c+it) \xi^{-c-it} \, dt
\end{equation}
is independent of $c\in (a,b)$ and satisfies $M M_c^{-1}g =g$.

For $\eps>0$ fixed, for
$f\in \cM^2_{(\eps)}\cap \cM^1_{(\eps)}$,
$f(x+iy)\to 0$ as $|y|\to+\infty$ uniformly in $x\in
(a,b)$.  
Therefore, 
$M_{\frac{n}{2}}^{-1} f$ is independent of $n=0,1,2,\dots$ and satisfies $M
M_{\frac{n}{2}}^{-1} f=f$.
We set
\begin{equation}\label{def-cM-inv}
M^{-1} f = M_{\frac{n}{2}}^{-1} f\,.
\end{equation}

Having constructed an inverse of $M$ on a dense subspace of
$\cM^2$, if we show that on this subspace
\begin{equation}\label{second-isometry}
\|M^{-1} f \|_{L^2( (0,+\infty),\,\frac{e^{2\xi}}{\xi}\,d\xi)} = \| f\|_{\cM^2} \,,
\end{equation}
the conclusion will follow.

For $f\in  L^2\big(\{c\}+i\bR\big)$ define $\vp_c(\xi)=
\xi^{-c} f\big( c-i\log \xi)$ .  Then, setting $L^2_c = L^2\big(
(0,+\infty), t^{2c}\frac{dt}{t}\big)$, 
\begin{equation}\label{phi-c}
\|\vp_c\|_{L^2_c} =\| f\|_{L^2(\{c\}+i\bR)}\,,
\end{equation}
since
\begin{align*}
\| \vp_c \|_{L^2_c}^2 
& = \int_0^{+\infty} |f\big( c-i\log \xi)|^2 \, \frac{d\xi}{\xi} 
= \int_{-\infty}^{+\infty} |f\big( c-it)|^2 \, dt 
 = \|f\|_{L^2 (\{c\}+i\bR)}^2\,.
\end{align*}

We claim that for such $f$  and $\xi>0$ we have
\begin{equation}\label{M-inv-c-f-equals-M-phi-c}
\xi^c M_c^{-1} f(\xi) =M \vp_c (c+i\log \xi)\,.
\end{equation}
For, 
\begin{align*}
M_c^{-1} f(\xi)
& = \frac{1}{\sqrt{2\pi}} \int_{-\infty}^{+\infty} f(c+it) \xi^{-c-it} \, dt\\
& = \frac{\xi^{-c} }{\sqrt{2\pi}} \int_0^{+\infty} f (c-i\log s) \xi^{i\log s}
\,\frac{ds}{s}\\
& = \frac{1}{\sqrt{2\pi}}  \int_0^{+\infty} s^c \vp_c (s) \xi^{i\log s}
\,\frac{ds}{s}\\
& = \xi^{-c}  M \vp_c (c+i\log \xi)\,,
\end{align*}
as we claimed.
\ms

Lemma 2.3 in \cite{BuJa}  shows that
$$
M: L^2_c \to  L^2\big(\{c\}+i\bR\big)
$$
is an isometry.  Finally, using \eqref{phi-c} and \eqref{M-inv-c-f-equals-M-phi-c} we have
\begin{align*}
\| M^{-1} f\|_{L^2( (0,+\infty),\,\frac{e^{2\xi}}{\xi}\,d\xi)}^2
& = \sum_{n=0}^{+\infty} \frac{2^n}{n!} 
\int_0^{+\infty} |M^{-1} f(\xi)|^2 \xi^n \, \frac{d\xi}{\xi}
\, = \, \sum_{n=0}^{+\infty} \frac{2^n}{n!} 
\int_0^{+\infty} |M_{\frac{n}{2}}^{-1} f(\xi)|^2 \xi^n \, \frac{d\xi}{\xi}\\
& = \sum_{n=0}^{+\infty} \frac{2^n}{n!} 
\int_0^{+\infty} \big|M \vp_{\frac{n}{2}} \big(\tnt+i\log\xi\big)\big|^2 \,
\frac{d\xi}{\xi} 
\displaystyle{\, 
=\,  \sum_{n=0}^{+\infty} \frac{2^n}{n!} 
\int_{-\infty}^{+\infty} \big|M \vp_{\frac{n}{2}} (\tnt+iy)\big|^2 \,
dy }\\
& = \sum_{n=0}^{+\infty} \frac{2^n}{n!} 
\int_0^{+\infty} |\vp_{\frac{n}{2}} (\xi)|^2 \xi^n \, \frac{d\xi}{\xi}
\, =\,  \sum_{n=0}^{+\infty} \frac{2^n}{n!} 
\int_{-\infty} ^{+\infty} |f (\tnt+iy)|^2 \, dy\\
& = \|f\|_{\cM^2}^2\,.
\end{align*}
This completes the proof.
\qed
\ms

We conclude this section discussing some property of the space
$\cH:= \frac{2^z}{\Gamma(1+z)} \cM^2$.

\proof[Proof of Theorem \ref{new}.]
We have already proved that $M_\Delta: A^2(\Delta)\to\cH$ is a
surjective isometry.  Proposition \ref{Mf-exp-type} shows that
elements of $\cH$ are holomorphic functions in 
$\cR$ of exponential type at most $\pi/2$, while Example
\ref{imag-polys} easily implies that the polynomials are dense in
$\cH$. 
  
Clearly, the reproducing kernel of $\cH$ is given by
\begin{align*}
H(z,w) & = \frac{2^z}{\Gamma(1+z)} 
\Big(\frac{1}{2\pi}\frac{\Gamma(z+\ov w)}{2^{z+\ov w}}\Big)
\frac{2^{\ov w}}{\Gamma(1+\ov w)}  = 
\frac{1}{2\pi} \frac{\Gamma(z+\ov w)}{\Gamma(1+z)\Gamma(1+\ov w)} \,
. \qed
\end{align*}
\ms

\begin{remark}\label{rem-new}{\rm 
It is not difficult to show (see \cite{Bologna}) that $\cH$ contains
functions that are of exponential type $\pi/2-\eps$ for any $\eps>0$,  
and whose restriction to
the imaginary axis is again of exponential type.  This shows that no
Paley--Wiener type theorem can hold for $\cH$.  Moreover, the measure
on $\ov\cR$ that appears in the norm of $\cH$
$$
d\mu(x+iy):= 
\sum_{n=0}^{+\infty} \delta_{\frac{n}{2}}(x)\otimes
\frac{|\Gamma(\tnt+iy)|^2}{\Gamma(n+1)} 
dy
$$
is not translation invariant, in contrast to the invariance of
$\omega$.

The results of the next section concerning the zero-sets indicate that
it is easier to exploit the properties of $\cM^2$  than the ones of
$\cH$.  
}
\end{remark}
\ms

\section{Zero-sets 
}\label{zero-sets}
\ms

Denote by $\cE_\tau(\ov\cR)$ and $\cE_{<\tau}(\ov\cR)$ respectively, the space of
holomorphic functions on $\ov\cR$ that are of exponential type $\tau$
and of  exponential type less than $\tau$, respectively.
Let $\cZ(\cK)$ denote the collection of zero-sets for the space
$\cK$.  
\begin{prop}\label{zero-sets-containtment}
We have the inclusions
$$
\cZ(\cE_{<\frac\pi2}(\ov\cR))\subseteq
\cZ(\cM^2(\ov\cR )\subseteq 
\cZ(\cE_{\frac\pi2}(\ov\cR))\,.
$$
\end{prop}
\proof
From Proposition \ref{Mf-exp-type} 
and Theorem \ref{new} we have 
$$
\cM^2(\cR)\cap\Hol(\ov\cR)\subseteq \Gamma(1+\cdot) \cE_{\frac\pi2}(\ov\cR) \,.
$$
The conclusion now follows from the above inclusion and Proposition
\ref{prop-1.2} (iii). \ms
\qed

In order to prove Theorem \ref{zero-set-thm} we need a couple of preliminary results
that may be of independent interest.
Recall that  the elementary
Weierstrass factor $E(z,p)$ equals $1-z$ when $p=0$, 
while  $E(z,p)= (1-z)e^{z+\cdots+\frac{z^p}{p}}$,
for a positive integer $p$.

\begin{prop}\label{type-can-pro-prop}
Let $\{z_j\}_{j=1,2\dots}\subseteq \cR$, with $|z_j|\to+\infty$,
$|z_j|>1$ 
 having exponent of convergence $1$ and assume $d^+<\infty$. 
Set $z_{-j}=-z_j$ and consider the 
sequence $\{z_j\}_{j\neq0}$. 
Then the infinite product
$\prod_{j\neq0}  E(z/z_j,1)$ converges to an entire function $\Pi(z)$
of exponential type at most $\pi d^+$.
\end{prop}

\proof
The sequence $\{z_j\}_{j\neq0}$ has exponent of convergence $1$
so that the product $\prod_{j\neq0}  E(z/z_j,1)$ converges to an
entire function $\Pi(z)$.  Thus, we only need to prove the statement
about the exponential type.

Since
$E(z/z_j,1)E(-z/z_j,1)=E(z^2/z_j^2,0)$ we have
\begin{align*}
\log \big| \prod_{j\neq0}  E(z/z_j,1)\big|
& = \log \big| \prod_{j=1}^{+\infty}  E(z^2/z_j^2,0)\big| \\
& = \sum_{j=1}^{+\infty} \log \big|E(z^2/z_j^2,0)\big| \\
& \le \sum_{j=1}^{+\infty} \log \big( 1+ |z^2/z_j^2| \big) \,.
\end{align*}
Setting $|z|=r$ we then have
\begin{align*}
\log |\Pi(z)| 
& \le  \sum_{j=1}^{+\infty} \log \big( 1+ r^2/|z_j|^2| \big) \\
& = \int_1^{+\infty} \log \big( 1+ r^2/t^2| \big)\, dn(t) \\
& = \log \big( 1+ r^2/t^2 \big) n(t)\Big|_1^{+\infty} + 2r^2 
\int_1^{+\infty} \frac{n(t)}{t^3(1+ r^2/t^2)} \, dt \,.
\end{align*}
Notice that
\begin{align*}
\log \big( 1+ r^2/t^2 \big) n(t)\Big|_1^{+\infty} 
 = \lim_{R\to+\infty} \log \big( 1+ r^2/R^2 \big) n(R)  
 \le C \lim_{R\to+\infty} (r^2/R^2)R = 0\,.
\end{align*}
Therefore,  given $\eps>0$ there exists $A>0$ large enough so that
\begin{align*}
\log |\Pi(z)| 
& \le 2r^2 
\int_1^{+\infty} \frac{n(t)}{t(t^2+ r^2)} \, dt \,=\,  2r^2 \bigg( \int_1^A + \int_A^{+\infty} \bigg) \frac{n(t)}{t(t^2+
  r^2)} \, dt \\
& \le  2r^2 \bigg(  C \int_1^A \frac{1}{t^2+ r^2} \, dt +
(d^++\eps) \int_A^{+\infty} \frac{1}{t^2+ r^2} \, dt \bigg) \\
& \le  2r^2\bigg(  C \frac{A}{1+r^2} \
+ \frac{(d^++\eps)}{r}  \int_0^{+\infty} \frac{1}{1+s^2}\, ds \bigg) \\
& \le CA +\pi(d^++\eps)r\,,
\end{align*}
as $r\to+\infty$, and $\Pi(z)$ is of exponential type at most $\pi d^+$.
\qed\ms

In order to obtain a necessary condition for a sequence to be a
zero-set for $\cM^2$, we recall the classical Carleman formula
for the right half-plane.

\begin{thm}\label{classical-Carleman} {\rm {\bf (Carleman)}}
Let $f\in \Hol(\ov\cR)$ and let $\{z_j\}$ its zero-set,
with $r_j\ge1$, where $z_j=r_j e^{i\theta_j}$.
Then, for $R\ge1$ we have
\begin{align}
& \sum_{r_j\le R} \Big(\frac{1}{r_j}-\frac{r_j}{R^2}\Big)\cos \theta_j
\notag \\
& \qquad = \frac{1}{2\pi} \int_1^R \Big(\frac{1}{y^2}-\frac{1}{R^2}\Big)
\log\big| f(iy)f(-iy)\big|\, dy 
+ \frac{1}{\pi R} \int_{-\frac\pi2}^{\frac\pi2} \log | f(Re^{it})|\,
\cos t\, dt +A(R) \,, \label{Car-cla-for}
\end{align}
where $A(R)$ is a bounded function of $R$.
\end{thm}

\begin{prop}\label{our-Carleman}
Let $f\in \Hol(\ov\cR)$ be such that
\begin{itemize}
\item[(i)] $\sup_{-\frac\pi2\le\theta\le\frac\pi2} |f(Re^{i\theta})|
  \le c_1e^{BR\log R}$ with $B>0$ and some $c_1>0$;\smallskip
\item[(ii)] $|f(iy)|\le c_2 e^{A|y|}$ with $A>0$, and some $c_2>0$. 
\end{itemize}
If $\{z_j\}$ are the zeros of $f$ with $|z_j|\ge1$, then
\begin{equation}\label{our-Car-est}
\sup_{R>0} \frac{1}{\log R} \sum_{r_j\le R}
\Big(\frac{1}{r_j}-\frac{r_j}{R^2}\Big)\cos \theta_j \le \frac1\pi
(A+2B) \,.
\end{equation}
\end{prop}
\proof
Denote by $ I(R) +J(R)+A(R)$ the right hand side in
\eqref{Car-cla-for}.  Consider $I(R)$ and set
$$
I_\pm (R) = \frac{1}{2\pi} 
\int_1^R  \Big(\frac{1}{y^2}-\frac{1}{R^2}\Big) 
 \log^\pm\big| f(iy)f(-iy)\big|\, dy \,.
$$
Clearly we have
$I_- (R)\le 0\le I_+(R)$.  Moreover, since (ii)  implies $\log^+ 
|f(iy)f(-iy)|\le 2A|y|+ C$, it follows that
\begin{align*}
0\le \frac{1}{\log R} I_+(R)
& \le  \frac{1}{2\pi\log R} \int_1^R
\Big(\frac{1}{y^2}-\frac{1}{R^2}\Big)(Ay+C) \, dy 
\,,
\end{align*}
which tends to $\frac{A}{\pi}$, as $R\to+\infty$.

Next we consider $J(R)/\log R$.  We observe that since $f$ vanishes of
finite order in $\cR$ we have
$$
\Big| \int_{-\frac\pi2}^{\frac\pi2} \log_- |f(Re^{i\theta})|
\cos\theta\, d\theta \Big|\le C\,.
$$
Therefore, using (i) we have
\begin{align*}
\bigg| \frac{J(R)}{\log R}\bigg|
& \le \frac{C}{\pi R \log R} 
+ \frac{1}{\pi R \log R} \int_{-\frac\pi2}^{\frac\pi2} \log_+ |f(Re^{i\theta})|
\cos\theta\, d\theta \\
& \le \frac{C}{\pi R \log R} + \frac{2B}{\pi}\,,
\end{align*}
which tends to $\frac{2B}{\pi}$, as $R\to+\infty$.

Now notice that the left-hand side of \eqref{Car-cla-for}  is
non-negative and increasing in $R$.  Therefore,
$$
0\le \limsup_{R\to+\infty} \frac{1}{\log R} \Big(
I_+(R)+I_- (R) +J(R)+A(R) \Big) =
\frac{A}{\pi}
+ \frac{2B}{\pi} + \limsup_{R\to+\infty} \frac{I_- (R)}{\log R} \,.
$$ 
Thus, $I_- (R)/\log R$ must remain bounded from below it follows that
\begin{equation}\label{our-Car-est-limsup}
\limsup_{R\to+\infty} \frac{1}{\log R} \sum_{r_j\le R}
\Big(\frac{1}{r_j}-\frac{r_j}{R^2}\Big)\cos \theta_j \le \frac1\pi
(A+2B) \,. \qed
\end{equation}
\ms

\proof[Proof of Theorem \ref{zero-set-thm}]

(i) Using Proposition 
\ref{type-can-pro-prop}
we can construct an entire function $\Pi$ of exponential type
$\tau<\frac\pi2$, whose zeros in $\ov\cR$ is exactly the sequence
$\{z_j\}$.  From Proposition \ref{prop-1.2} (iii), there exists $\delta>0$
such that $\Pi(z) \Gamma(1+\delta z)$ is in $\cM^2(\cR)\cap\Hol(\ov\cR)$.
\ms

(ii) 
From Theorem \ref{new} it follows that any $f\in \cM^2(\cR)\cap\Hol(\ov\cR)$ 
 satisfies
the hypotheses of Proposition \ref{our-Carleman},  
so that using the asymptotics for the Gamma function
\eqref{asy-basic}, 
 conclusion \eqref{our-Car-est} holds for $f$ with $A=\eps$, for every $\eps>0$
and $B=1$. 

Observe that, for $0\le\delta<1$, 
$$
\frac{1}{r_j} -\frac{r_j}{R^2} \ge \frac{\delta}{r_j}\quad
\text{if and only if}\quad
r_j\le R\sqrt{1-\delta}\,.
$$
Therefore,
\begin{align*}
\sum_{r_j\le R} \Big(\frac{1}{r_j}-\frac{r_j}{R^2}\Big)\cos \theta_j
& \ge \sum_{r_j\le R\sqrt{1-\delta}} \Big(\frac{1}{r_j}-\frac{r_j}{R^2}\Big)\cos \theta_j
\\
& \ge \delta \sum_{r_j\le R\sqrt{1-\delta}} \frac{1}{r_j} \cos \theta_j \,,
\end{align*}
so that,
\begin{align*}
\limsup_{R\to+\infty} \frac{1}{\log R} \sum_{r_j\le R} \Big(\frac{1}{r_j}-\frac{r_j}{R^2}\Big)\cos \theta_j
& \ge \delta  \limsup_{R\to+\infty} \frac{1}{\log R}
\sum_{r_j\le R\sqrt{1-\delta}} \frac{1}{r_j} \cos \theta_j \\
& = \delta  \limsup_{R'\to+\infty} \frac{1}{\log R'- \log\sqrt{1-\delta}}
 \sum_{r_j\le R'} \frac{1}{r_j} \cos \theta_j \\
& = \delta  \limsup_{R\to+\infty} \frac{1}{\log R}
 \sum_{r_j\le R} \frac{1}{r_j} \cos \theta_j 
\,.
\end{align*}
The conclusion now follows.
\qed
\ms

\begin{remark}{\rm 
Since the Hardy space $H^2(\cR)$ is contained in $\cM^2(\cR)$, so that
$\cZ(H^2(\cR))\subseteq \cZ( \cM^2(\cR))$ it follows that if 
$\{z_j\}$ is such that $|z_j|\to+\infty$ and 
$$
\sum_j
\frac{\Re z_j}{1+|z_j|^2} <\infty\, ,
$$
then $\{z_j\}$ is also a zero-set for $\cM^2(\cR)\cap\Hol(\ov\cR)$.  
We now show that:\smallskip
\begin{itemize}
\item[(a)] $\cZ(H^2(\cR))\subsetneq \cZ( \cM^2(\cR))$; \smallskip
\item[(b)] there exist sequences in $\cZ( \cM^2(\cR))$ (actually, in
  $\cZ(H^2(\cR))$) that do not satisfy condition (i) in Theorem
  \ref{zero-set-thm}; hence, this condition is not
necessary. \ms
\end{itemize}

\noindent 
(a) Any sequence  $\{z_j\}$ satisfying condition (i) in
Theorem \ref{zero-set-thm}, 
contained in a sector  $|\arg z_j|\leq \vartheta<\frac \pi 2$ and such that $\sum_j
\frac{1}{|z_j|}=+\infty$, is in $\cZ( \cM^2(\cR))$ but not in
$\cZ(H^2(\cR))$.  As special cases, we find all sequences of the form
$z_j = aje^{i\vartheta}$, $j=1,2,\dots$, with $a>2$.  \ms

\noindent 
(b)  Any sequence  $\{z_j\}$ contained in a strip $S_b$ with
exponent of convergence $\rho>1$ and such that  $\sum_j
\Re(1/z_j)<\infty$, is in $\cZ(H^2(\cR))$ hence in  $\cZ(
\cM^2(\cR))$ but does not satisfy condition (i) in
Theorem \ref{zero-set-thm}.  As special cases, we find all sequences of the form
$z_j = a+i j^\alpha$, $j=1,2,\dots$, with $\frac12 <\alpha<1$.  \qed

\ms

\ms
}
\end{remark}

\section{Sets of uniqueness and 
the  \MS\ problem for the Bergman space}\label{MS-sec}\ms 

The \MS\ problem for the Bergman
space was formulated by S. Krantz, C. Stoppato and the first author, see 
\cite{KPS2}.
 In Theorem 3.1 of the same paper, using an {\em ad hoc} method,  
it is shown that if $\lambda_k = \eps_0 +ak +ib$,
 where $\eps_0>0$, $b\in\bR$ and $0<a<1$, the set
 $\{\z^{\lambda_k-1}\}$ is a complete set in $A^2(\Delta)$.

We recall a classical result by Fuchs \cite{Fuchs}, concerning sets
of uniqueness for functions that are of exponential type in a
half-plane.  It says that 
if $\lambda_k>0$, $\lambda_{k+1}-\lambda_k>\delta>0$ and the 
sequence $\{\lambda_k\}$ has lower density $d^->\frac12$, then 
 $\{\lambda_k\}$ is a set of uniqueness for the functions that are of
 exponential type $\frac\pi2$ in $\cR$.   Since the sets of uniqueness
 satisfy the inverse inclusions of Proposition
 \ref{zero-sets-containtment},  it follows that
the set 
$\{\z^{\lambda_k-1}\}$ is a complete set in $A^2(\Delta)$. \ms

In light of Proposition
 \ref{zero-sets-containtment}, it is clear that Theorem 3.1 in
 \cite{KPS2} follows from Fuchs' result.
Our Theorem \ref{MS-thm} applied to a sequence $\{\lambda_k\}$ 
on the positive half-line
gives that if it has lower density greater than $\frac2\pi$,
regardless of being separated or not, then $\{\lambda_k\}$ is a set of
uniqueness for $\cM^2$.  Thus, while it does not contain Fuchs' result,
it is of much more general nature, since it only assumes 
a lower bound of a quantity depending on the number of the
$\{\lambda_k\}$ in the half-disks $\{|z|\le R, \Re z>0\}$ and not on
any type of distribution of the $\lambda_k$'s. \ms

\proof[Proof of Theorem \ref{MS-thm}]
Let $\{z_j\}\subseteq \{ \Re z \ge \eps_0\}$ be such that
\eqref{our-Carleman-cond} is violated and let $f\in\cM^2$ vanish on $\{z_j\}$.  
If $f\in \Hol(\ov\cR)$ the result follows at once from Theorem
\ref{zero-set-thm} (ii).   

For a generic $f\in\cM^2$ and $0<\eps<\eps_0$,
 we consider the points $w_j=z_j-\eps$ 
 and the function $f_\eps(z) =f(z+\eps)$.  Then $f_\eps
\in \cM^2\cap\Hol(\ov\cR)$ and vanishes at the $w_j$'s.  If we show
that the sequence $\{w_j\}$ violates condition
\eqref{our-Carleman-cond}, it would follow that $f_\eps$, hence $f$,
is identically zero.  Notice that $|w_j|=|z_j-\eps|<|z_j|$ so that
$$
\Re \rw  \ge \frac{ \Re z_j -\eps}{|z_j|^2} 
= \Re \rz
- \frac{\eps}{|z_j|^2}\, .
$$
Therefore,
\begin{align*}
\sum_{|z_j|\le R} \Re \rz 
& \le \sum_{|z_j|\le R}  \Re \rw 
+ \frac{\eps}{|z_j|^2}\\
& \le \sum_{|w_j|\le R}  \Re \rw 
+ \sum_{|z_j|\le R}  \frac{\eps}{|z_j|^2} \\
&  \le \sum_{|w_j|\le R}  \Re \rw 
+ \frac{\eps}{\eps_0} \sum_{|z_j|\le R}  \Re \rz 
\, .
\end{align*}
Hence, 
$$
\Big(1-\frac{\eps}{\eps_0}\Big) 
\limsup_{R\to+\infty} \frac{1}{\log R} 
\sum_{|z_j|\le R} \Re \rz 
\le \limsup_{R\to+\infty} 
\frac{1}{\log R}\sum_{|w_j|\le R}  \Re \rw \, ,
$$
for every $0<\eps<\eps_0$.
The conclusion now follows by taking $\eps>0$ sufficiently small.
\qed
\ms

\ms

\section{Proof of Theorems \ref{uniqueness-measure}
and \ref{proj-thm}.}\label{proj-thm-sec}\ms 

We remark that in the next proof in particular we show how the measure
$\omega$ was determined.
 
\proof[Proof of Theorem  \ref{uniqueness-measure}.]
By Theorems \ref{reduction-thm} and \ref{mapping-M} it suffices to 
show that, if the Mellin transform $M: L^2 \big((0,+\infty),
\frac{e^{2\xi}}{\xi} d\xi \big) \to L^2(\ov\cR,d\mu)$ is an isometry
and $\mu$ is translation invariant, then $\mu=\omega$.

Assume then that $d\mu=d\nu(x)\otimes dy$.  Let 
$\eta_{1,\eps},\eta_{2,\eps'}$ be
smooth cut-off functions with support in 
$[\eps,1/\eps]$ and 
$[-1/\eps',1/\eps']$, respectively, and identically 1 in   
$[2\eps,1/(2\eps)]$ and 
$[-1/(2\eps'),1/(2\eps')]$, respectively.    Then
$$
\eta_{1,\eps}(x)\eta_{2,\eps'}(y) 
\, d\nu(x)\otimes dy \to d\mu(x,y)
$$
as Borel measures, as $\eps,\eps'\to0^+$. 

Notice that for $\vp \in  C_0^\infty(0,+\infty)$, $M\vp=\cF^{-1}
(\vp\circ\exp)$ is an entire 
function.
Then,  
for any 
$\vp,\psi\in C_0^\infty(0,+\infty)$  we have
\begin{align}
&\la \vp,\psi \ra_{L^2((0,+\infty),\frac{e^{2\xi}}{\xi}d\xi)} \notag \\
& = \iint_{\ov\cR} M\vp(x+iy)\ov{M\psi(x+iy)}\,  d\nu(x) \, dy  \notag \\
& = \lim_{\eps,\eps'\to0^+} \iint_{\ov\cR} M\vp(x+iy)\ov{M\psi(x+iy)}\,
\eta_{1,\eps}(x)\,  d\nu(x) \, \eta_{2,\eps'}(y)dy \notag \\
& =   \frac{ 1}{2\pi}    \lim_{\eps,\eps'\to0^+} \int_0^{+\infty}\int_0^{+\infty} \vp(s)\ov{\psi(t)} 
\Big( \iint_{\ov\cR} s^{x+iy} t^{x-iy} \,
\eta_{1,\eps}(x)\,  d\nu(x)\, \eta_{2,\eps'}(y)dy \Big)\, \frac{ds}{s}
\frac{dt}{t} \, , \label{ci-serve?}
\end{align}
where we can interchange the integration order since the integrals
converge absolutely.  The inner integral in the right hand side above
equals, using Fubini's theorem again,
\begin{align*}
&   \frac{ 1}{2\pi}  \int_0^{+\infty} \int_{-\infty}^{+\infty} s^{iy} t^{-iy}
\eta_{2,\eps'}(y)
s^x t^x \eta_{1,\eps}(x)\, dy \,  d\nu(x)\\
& =   \frac{ 1}{\sqrt{2\pi}} \cF(\eta_{2,\eps;})\big(\log(s/t)\big)  \int_0^{+\infty} e^{x\log(st)} 
\eta_{1,\eps}(x)\,  d\nu(x) \\
& =  \cF(\eta_{2,\eps'})\big(\log(s/t)\big) 
\cF\big( \eta_{1,\eps}\nu\big) \big(i\log(st)\big) 
\, ,
\end{align*}
 observing that $\cF\big( \eta_{1,\eps}\nu\big)$ is
entire.  
Plugging  this equality into the right hand side of \eqref{ci-serve?}
we find that
\begin{align}
&\la \vp,\psi \ra_{L^2((0,+\infty),\frac{e^{2\xi}}{\xi}d\xi)} \notag \\
& = \lim_{\eps,\eps'\to0^+} \int_0^{+\infty}\int_0^{+\infty} \vp(s)\ov{\psi(t)} 
\cF(\eta_{2,\eps'})\big(\log(s/t)\big) 
\cF\big( \eta_{1,\eps}\nu\big) \big(i\log(st)\big) \, \frac{ds}{s}
\frac{dt}{t} \notag \\
& = \lim_{\eps,\eps'\to0^+} \int_{-\infty}^{+\infty}\int_{-\infty}^{+\infty}
\vp(e^y)\ov{\psi(e^v)} \cF(\eta_{2,\eps'}) (y-v) 
\cF\big( \eta_{1,\eps}\nu\big) (i(y+v))  \,dydv \notag \\
& = \lim_{\eps\to0^+}  
\sqrt{2\pi} \int_{-\infty}^{+\infty}
\vp(e^y)\ov{\psi(e^y)} 
\cF\big( \eta_{1,\eps}\nu\big) (2iy)  \,dy\notag \\
& = \lim_{\eps\to0^+}  \sqrt{2\pi}  \int_0^{+\infty} \vp(\xi) \ov{\psi(\xi)} \cF\big( \eta_{1,\eps}\nu\big) ( 2i\log \xi)
\frac{d\xi}{\xi}\, .
\end{align} 

Hence, if such a measure $\nu$  exists, for $\xi>0$ it
must satisfy the identity 
\begin{align}
 e^{2\xi} 
& =\sqrt{2\pi} \lim_{\eps\to0^+} \cF\big(
\eta_{1,\eps}\nu\big) ( 2i\log \xi) \notag\\
& = 
\lim_{\eps\to0^+} \int_0^{+\infty} e^{2x\log \xi}
\eta_{1,\eps} (x) 
\, d\nu(x)   \notag \\
& =  \int_0^{+\infty} e^{2x\log \xi}
\, d\nu(x)  \label{new-tag}
\end{align}
by the monotone convergence theorem.  This in particular implies that: 
\begin{itemize}
\item[{\tiny $\bullet$}]
$\nu(\{0\})=1$ (by letting $\xi\to 0^+$) and
$\nu\big([0,+\infty)\big)=e^2$; \smallskip 
\item[{\tiny $\bullet$}]
$x\mapsto \xi^{2x}\in L^1(d\nu)$, for all $\xi>0$; \smallskip 
 \item[{\tiny $\bullet$}] for all $s\in\bR$,
$$
\exp\{2e^{\frac s2}\}
=  \int_0^{+\infty} e^{xs}
\, d\nu(x)  \, .
$$
 \end{itemize} 

Hence, $\cF \nu$ can be extended to the entire function
$\cF \nu(t+is) =  \frac{1}{\sqrt{2\pi}}\exp\{2e^{-i\frac{t+is}{2}} \}$,
so that 
\begin{equation}\label{last-claim}
\nu =
\frac{1}{\sqrt{2\pi}} \cF^{-1} \big(  \exp\{2e^{-\frac{i}{2}(\cdot)} \}\big)\, .
\end{equation}
 
Thus, if $\nu$ exists, it is unique and it is given by
\eqref{last-claim}.  

Observe that the function $\exp\{2e^{-\frac{i}{2}x}\}$ is bounded on
$\bR$, so that its (inverse) Fourier transform is well defined as a tempered
distribution.  
Moreover, it is the restriction to the real line of the
entire function
$$
G(z) = \exp\{2e^{-\frac{i}{2}z}\} = \sum_{n=0}^{+\infty}
\frac{2^n}{n!} e^{-i\frac{n}{2}z}\, ,
$$ 
and by \cite{Lukacs} Theorem 2 we know that $\cF^{-1} G$ is a finite Borel
measure with support in $[0,+\infty)$,  that in fact we can compute
explicitly.   
The series $\sum_{n=0}^{+\infty}
\frac{2^n}{n!} e^{-i\frac{n}{2}x}$ converges uniformly on compact
subsets of the line to the bounded function
$\exp\{2e^{-\frac{i}{2}x}\}$, hence in the sense of tempered distributions.
Therefore, 
\begin{align*}
\nu (x) & = \frac{1}{\sqrt{2\pi} } \cF^{-1} \Big( \sum_{n=0}^{+\infty}
\frac{2^n}{n!} e^{-i\frac{n}{2}x} \Big) (x) = \sum_{n=0}^{+\infty}
\frac{2^n}{n!} \delta_{\frac{n}{2}}(x)\, ,
\end{align*}
as we wished to show.
\qed\ms

\proof[Proof of Theorem  \ref{proj-thm}.]
For $1<p<\infty$, the dual of $\cM^p$ with respect to the
$L^2(\ov\cR,d\omega)$-inner product is $\cM^{p'}$, with $1/p+1/p'=1$.
If the projection $P:L^p(\ov\cR,d\omega)\to \cM^p$ were bounded, $K_w
\in \cM^{p'}$, for any $w\in\cR$.  Thus, since $P$ is self-adjoint, it
suffices to show that $K_w\not\in L^p(\ov\cR,d\omega)$ for any
$p>2$. 

This is a simple  application of the asymptotics of the Gamma function.
For $w=u+iv\in\cR$ fixed
\begin{align*}
\|K_w\|_{\cM^p}^p 
& = \frac{1}{(2\pi)^p}  \Big\| \frac{\Gamma(\cdot+\ov w)}{ 2^{(\cdot+\ov w)}} \Big\|_{\cM^p}^p
\\
& = \frac{1}{(2\pi)^p 2^{pu}} \sum_{n=0}^{+\infty} \frac{2^{n(1-\frac{p}{2})}}{n!} \int_{-\infty}^{+\infty} 
\big| \Gamma(\textstyle{\tnt+u+iy})\big|^p \, dy \,.
\end{align*}
Now, using \eqref{asy-basic} we see that there exists an absolute
constant $C>0$ such that
\begin{align*}
\big| \Gamma(\textstyle{\tnt+u+iy})\big|^p 
& \ge C \exp\Big\{  p\Big[ \textstyle{ \big( \frac{n-1}{2}+u\big)
  \log\big( (\frac{n}{2} +u)^2 +y^2 \big)^{1/2}  -\frac{n}{2} -u -
  |y|\arctan (\frac{|y|}{\frac{n}{2}+u}) }
  \Big] \Big\} \\
& \ge  Ce^{-pu}  
\exp\Big\{  p\Big[ \textstyle{ \big( \frac{n-1}{2}+u\big)
  \log (\frac{n}{2} +u)  -\frac{n}{2}  -
  \frac{\pi}{2}|y| }
  \Big] \Big\} \\
& \ge Ce^{-pu}  e^{-p\frac{\pi}{2}|y|}  
\exp\Big\{  \textstyle{ \frac{p(n-1)}{2}\log (\frac{n}{2})
  -\frac{pn}{2} } \Big\} \, .
\end{align*}
Therefore,
\begin{align*}
\|K_w\|_{\cM^p}^p 
& \ge C_u \sum_{n=0}^{+\infty} \frac{2^{n(1-\frac{p}{2})}}{n!} 
\exp\Big\{  \textstyle{ \frac{p(n-1)}{2}\log (\frac{n}{2})
  -\frac{pn}{2} } \Big\} \\ 
& \ge C_u \sum_{n=0}^{+\infty} \frac{2^{n(1-p-\frac{p}{2}\log
    2)}}{n!\,n^{\frac{p}{2}}} e^{\frac{p}{2}n\log n}  \, ,
\end{align*}
which clearly diverges when $p>2$. \qed \ms

\bibliography{m2-mathscinet-1}
\bibliographystyle{abbrv}

\end{document}